\documentclass[a4paper,12pt]{amsart}
\usepackage[letterpaper, margin=1.2in]{geometry}
\usepackage{latexsym}
\usepackage{amssymb}
\usepackage{amsmath}
\usepackage{amsfonts}
\usepackage{color}
\usepackage{verbatim}
\everymath{\displaystyle}
\newtheorem{thm}{Theorem}[section]

\newtheorem{lem}[thm]{Lemma}

\newtheorem{prop}[thm]{Proposition}
\theoremstyle{definition}

\newtheorem{rem}{Remark}

\def\fph{\mathbb{F}_{\ph}}

\newcommand{\Z}{\mathbb Z}
\newcommand{\z}{\mathbb Z}
\newcommand{\Q}{\mathbb Q}

\newcommand{\F}{\mathbb F}

\def\F{\mathbb{F}}

\newcommand{\p}{\mathfrak{p}}
\def\ol{\overline}
\def\al{\alpha}
\def\la{\lambda}

\def\om{\omega}

\def\md#1{\ \mbox{\rm(mod }{#1})}
\def\nph#1{N_{\ph}(#1)}
\def\npp#1{N_{\ph}^-(#1)}
\def\ph{\phi}
\def\ind{\text{ind}}
\def\tr{\triangle}
\newcounter{cs}
\stepcounter{cs}
\newcommand{\casos}{\begin{itemize}}
\newcommand{\fcasos}{\end{itemize}\setcounter{cs}{1}}

\newfont{\tit}{cmr12 scaled \magstep3}
\begin{document}
\title[]{On  the monogenity of quartic number fields defined by  $x^4+ax^2+b$}
\author{Lhoussain El Fadil and Istv\'an Ga\'al}{}
\address{Faculty of Sciences Dhar El Mahraz, P.O. Box  1796 Atlas-Fez , Sidi mohamed ben Abdellah University,  Morocco}\email{lhoussain.elfadil@usmba.ac.ma}
\address{Institute of Mathematics, University of Debrecen, H-4002 Debercen, Pf.400, Hungary
}\email{gaal.istvan@unideb.hu}

\begin{abstract}
For any quartic number field $K$ generated by a root $\alpha$ of an irreducible trinomial 
of type $x^4+ax^2+b\in \Z[x]$, we  characterize when $\Z[\alpha]$ is integrally closed. 
Also for $p=2,3$, we explicitly give
the highest power of $p$ dividing $i(K)$, the common index divisor of $K$.
For a wide class of monogenic trinomials of this type we prove that 
up to equivalence there is only one generator of power integral bases in $K=\Q(\alpha)$.
We illustrate our statements with a series of examples.
\end{abstract}


\maketitle

\vspace{0.3cm}
\section{Introduction}
\subsection{Monogenity of number fields and polynomials}
Let $K$ be a number field of degree $n$ with ring of integers $\Z_K$, and absolute discriminant $d_K$.
The number field $K$ is called {\it monogenic} if it admits a
{\it power integral basis}, that is an integral basis of type $(1,\alpha,\ldots,\alpha^{n-1})$
for some $\alpha\in\Z_K$.
Monogenity of number fields is a classical problem of algebraic number theory, going back to Dedekind, Hasse and Hensel, cf e.g. \cite{Ha,  He} and \cite{G19} for the present state of this area. 
It is called a problem of Hasse to give an arithmetic characterization of those number fields which have a power integral basis \cite{Ha, He, MNS}.
For any primitive element $\alpha$ of $\Z_K$ (that is $\alpha\in\Z_K$ with $K=\Q(\alpha)$)
we denote by 
\[
\ind(\alpha)=(\Z_K:\Z[\alpha])
\]
the {\it index of $\alpha$}, that is 
the index of the $\Z$-module $\Z[\alpha]$ in the free $\Z$-module $\Z_K$ of rank $n$.
As it is known \cite{G19}, we have
\[
\triangle(\alpha)=\ind(\alpha)^2\cdot d_K
\]
where $\tr(\alpha)$ is the discriminant of $\alpha$. If $F(x)$ is the minimal polynomial of $\alpha$ then we have $\Delta(F)=\Delta(\alpha)$ and
we also use the {\it index of} $F$, $\ind(F)=\ind(\alpha)$.
We also say that the polynomial $F(x)\in\Z[x]$ is {\it monogenic}, if $\ind(F)=1$, that is a root $\alpha$ 
of $F(x)$ generates a power integral basis in $K=\Q(\alpha)$.
Obviously, if $F(x)$ is monogenic, then $K$ is also monogenic, but the converse is not true:
there may exist a power integral basis in $K$, even if $F(x)$ is not monogenic.
Further, note that $F(x)$ is monogenic (that is $(1,\alpha,\ldots,\alpha^{n-1})$
is an integral basis of $K$) if and only if $\Z[\alpha]$
is integrally closed.

The elements $\alpha$ and $\beta$ of $\Z_K$ are called {\it equivalent } if $\alpha\pm\beta\in\Z$.
Obviously, equivalent elements have the same indices.

Let $(1,\omega_1,\ldots,\omega_{n-1})$ be an integral basis of $K$.
The discriminant $\tr(L(x_1,\ldots,x_{n-1}))$ of the linear form 
$L(x_1,\ldots,x_{n-1})=\omega_1x_1+\ldots +\omega_{n-1}x_{n-1}$ can be written (cf. \cite{G19}) as
\[
\tr(L(x_1,\ldots,x_{n-1}))=(\ind(x_1,\ldots,x_{n-1}))^2\cdot d_K,
\]
where $\ind(x_1,\ldots,x_{n-1})$ is the {\it index form} corresponding to the 
integeral basis $(1,\omega_1,\ldots,\omega_{n-1})$, having the property that for any
$\alpha=x_0+\omega_1x_1+\ldots +\omega_{n-1}x_{n-1}\in\Z_K$ (with $x_0,x_1,\ldots,x_{n-1}\in\Z$)
we have $\ind(\alpha)=|\ind(x_1,\ldots,x_{n-1})|$.

Obviously, $\ind(\alpha)=1$ if and only if $(1,\alpha,\ldots,\alpha^{n-1})$ is an integral basis
of $K$. Therefore $\alpha$ is a {\it generator of a power integral bases} if and only if
$(x_1,\ldots,x_{n-1})\in\Z^{n-1}$ is a solution of the {\it index form equation}
\[
\ind(x_1,\ldots,x_{n-1})=\pm 1 \;\; {\rm in}\;\; (x_1,\ldots,x_{n-1})\in\Z^{n-1}.
\]

\subsection{Results on the monogenity of binomilas, trinomials}

The problem of testing the monogenity of number fields and constructing power integral bases have been intensively studied during the last decades, see \cite{G19} and the references therein.
An especially delicate and intensively studied problem is the monogenity of {\it pure fields}
$K$ generated by a root $\alpha$ of an irreducible polynomial $x^n-m$ (cf. e.g. \cite{AN, GR17, E6s, FG9}).

Recently, many authors are interested in the monogenity of {\it trinomials} and number fields defined by 
a root $\alpha$ of a three term irreducible polynomial 
$x^n+ax^m+b\in\Z[x]$.
Jakhar, Khanduja and Sangwan \cite{AK}, Jhorar and Khanduja \cite{SK} and
Jakhar and Kumar \cite{A6}
studied the integral closedness of $\Z[\al]$, $\al$ being a root of a trinomial. 
Their results were refined by Ibarra, Lembeck, Ozaslan, Smith and Stange \cite{Smt}.  
Recall that the results given in \cite{AK, SK}, can only decide if $\Z[\al]$ is integrally closed
($\al$ a root of the trinomial), but cannot decide whether the field 
$K=\Q(\al)$
is monogenic or not. 

Jones \cite{Lsc, Lac}, Jones and Tristan \cite{LPh} and Jones and White \cite{LD}
also investigated monogenity of some irreducible trinomials. 
In   \cite{Ga21} Ga\'al calculated all generators of power integral bases
of number fields defined by some sextic irreducible trinomials. 
Ben Yakkou and El Fadil \cite{BF} gave  sufficient conditions on coefficients of certain trinomials which guarantee the non-monogenity of the  number field defined by a root of such a trinomial. 
Also, El Fadil \cite{jnt6, E5}  gave necessary and sufficient conditions  on 
the monogenity of number fields, generated by certain quintic and sextic trinomials,
in terms of the coefficients $a$ and $b$ of the trinomials.

\subsection{The index of a number field}
The greatest common divisor of the indices of all integral primitive elements of $K$ is called the  index of $K$, and  denoted by $i(K)$.  A rational prime $p$ dividing $i(K)$ is called a {\it prime common  index divisor} of $K$. It is clear that if $\mathbb{Z}_K$ has a power integral basis, then the index of $K$ is trivial, namely $i(K)=1$. Therefore a field having a prime index divisor is not monogenic. 

 The first  number field with non trivial index was given by Dedekind in 1871,
 	 who exhibited examples in cubic and quartic number fields. For example, he considered the cubic field $K$ generated by a root of $x^3-x^2-2x-8$ and  showed that the prime $2$ splits completely in $\Z_K$. So, if $K$ 
were monogenic, then there would be a cubic polynomial, a root of which  generating $K$, that splits completely into distinct polynomials of degree $1$ in $\mathbb{F}_2[x]$. Since there are only $2$ distinct polynomials of degree $1$ in $\mathbb{F}_2[x]$, this is impossible. 

Based on these ideas and using Kronecker's theory of algebraic number fields, 
Hensel  gave necessary and sufficient conditions 
	for any prime integer  $p$ to be  a prime common index divisor \cite{He2}.
Hensel \cite{He} also showed that the prime divisors of $i(K)$ must be less than the
degree of the field $K$.	
	
For arbitrary number fields of degree $n\le 7$, Engstrom \cite{En2} characterized   $\nu_p(i(K))$,
the highest power of $p$ dividing $i(K)$,
 by the factorization of  $(p)$ into powers of prime ideals of $\Z_K$ for every positive prime $p\le n$. 
Problem 22 of Narkiewicz \cite{WN} asks for an explicit formula for the highest power of a
given prime integer $p$ dividing $i(K)$.

In \cite{Na}, Nakahara  studied the  index of non-cyclic but abelian biquadratic number fields. 
 In \cite{GPP}  Ga\'al, Peth\"o and Pohst characterized  the field indices of biquadratic number fields having Galois group $V_4$.  
In quintic number field $K$ defined by a trinomial $x^5 + ax^2 + b$, El Fadil \cite{E5} 
gave necessary and sufficient conditions for prime integer $p$ to divide $i(K)$
in terms of $a$, $b$.

\subsection{The purpose of the present paper}
In this paper we consider trinomials of type $x^4+ax^2+b$. This is a special type of trinomials,
but exactly these special properties enable us to formulate results that far exceed
some corresponding statements on general types of trinomials. 

Our paper was motivated by some results concerning trinomials of type $x^4 +ax+b$.
Alaca and Williams \cite{AW1}, \cite{AW2} constructed $p$-integral bases and integral bases
of quartic fields generated by a root of such a trinomial.
Davis and Spearman {\cite{DA}} characterized  the prime  divisors of 
quartic number fields $K$ generated by a root of 
such a trinomial.

 In this paper, for any  quartic number field $K$ generated by a root $\alpha$ of an irreducible 
trinomial $x^4 + ax^2 + b\in\Z[x]$, we give necessary and sufficient conditions 
for the integral closedness of $\Z[\alpha]$ in terms of the coefficients $a$ and $b$.
We also evaluate $\nu_p(i(K))$ for every prime integer $p$.
  Based on Engstrom's results given in \cite{En2},  the unique prime candidates to divide $i(K)$ are $2$ and $3$, that is $i(K)=2^{\nu_2}3^{\nu_3}$, with $0\le \nu_2\le 2$ and $0\le \nu_3\le 1$.  In \cite{DA}, for a quartic number field $K$ defined by a trinomial $x^4+ax+b\in \Z[x]$, and for every prime integer $p$, Davis and Spearman gave necessary and sufficient conditions on $a$ and $b$ so that $p$ is a common index divisor of $K$. Their method is based on the calculation of the $p$-index form of $K$, derived from $p$-integral bases of $K$.
Our results in Theorems \ref{npib2} and \ref{npib3}
are analogous to that given in \cite{DA}, but our method is totally different;
 the proofs of  Theorems \ref{npib2} and \ref{npib3} are based on prime ideal factorization, which is performed for quartic number fields in the thesis of Montes (1999). The first author is very thankful to Professor Enric Nart who provided him a copy of   Montes's thesis. 

\medskip

In the present  paper,  we consider three types of problems: For   any quartic number field $K$  generated by a root $\al$ of an irreducible trinomial $x^4+ax^2+b\in\Z[x]$,\\
--we characterize when   $\Z[\al]$  is integrally closed (Theorem \ref{intclos}),\\
--we give necessary and sufficient conditions for $2$ or $3$ to divide the index of $K$, in terms of $a$, $b$\\
--we determine $\nu_p(i(K))$ for $p=2,3$,\\
--we study  monogenity of $K$ (see Subsection $4.2$).\\
--for a wide class of monogenic trinomials of type $x^4+ax^2+b$ we show that up to equivalence
the root of the trinomial is the only generator of power integral bases.\\

We also provide a series of examples illustrating our results.
\\

	\section{Main Results}
Throughout this section   unless otherwise stated, $K$ is a number field generated by a root $\al$ of an
irreducible  trinomial $F(x)=x^4+ax^2+b\in \Z[x]$ and we assume 
that for every prime $p$,  $\nu_p(a)  <2$ or  $ \nu_p(b) < 4 $.

  Along this paper, for any integer $a\in \Z$ and a prime $p$
and we set  $a_p=\frac{a}{p^{\nu_p(a)}}$.
\\

\subsection{Integral closedness of $\Z[\al]$}
Our first theorem characterizes the integral closedness of $\Z[\al]$:
 \begin{thm}\label{intclos}
The ring  $\Z[\al]$ is the ring of integers of $K$  if and only if for every prime integer $p$, $p$ satisfies one of the  following conditions:
\begin{enumerate}
\item
 $\nu_p(a)\ge 1$ and  $\nu_p(b)=1$.
 \item
  $p=2$, $b\equiv 2\md4$ and  $a \equiv 3\md4$.
\item
   $p=2$  does not divide $b$ and $a \equiv 1-b\md4$.
   \item
   $p=2$  does not divide $ab$ and   $a \equiv 3\md4$ or $b\not\equiv a\md4$.
  \item
  $p\ge 3$, $p$ does not divide $b$ and  
	either $p^2$ does not divide $a^2-4b$ or $-a2^{-1}$ is not a square in  $\F_p$.
\end{enumerate}
   In particular, if for every prime integer $p$, $p$ satisfies one of these conditions, then $i(K)$ is trivial, that is $i(K)=1$.
 \end{thm} 

\smallskip

\subsection{The divisibility of the field index by $2$ and $3$}
Next, for $p=2,3$, we give necessary and sufficient  conditions 
for $p$ to divide the index of $K$ in terms of $a$, $b$.
Furthermore, in every case we explicitly give $\nu_p(i(K))$.
\begin{thm} \label{npib2}
The following table provides the value of  $\nu_2(i(K))$:

$$\begin{array}{|c|c|c|c|c|}
	\hline
 {\mbox conditions}& &&&i(K)\\
 \hline
	 \nu_2(a)\ge 1& \nu_2(a)=1&\nu_2(b)=3& b\equiv 48-4a\md{64}&1\\
       \cline{3-5} \nu_2(b)\ge 1&& \nu_2(b)=2k+1 (k\ge 2)& b\equiv -2^{2k}a\md{2^{2k+4}}&1\\
       \hline
    b\equiv 0\md2&   \nu_2(b) \mbox{odd}&a\equiv 7-b\md{8}&&1\\
     \cline{2-5}  a\equiv 1\md2 &&a\equiv 1\md4&b\equiv -4a\md{32}&1\\
	   \cline{3-5}   &  \nu_2(b)=2&a\equiv 3\md8&b\equiv 8-4a\md{16}&1\\
     \cline{3-5}   & &a\equiv 7\md8&b\equiv 16-4a\md{32}&1\\
     \cline{3-5}   & &a\equiv 7\md8&b\equiv -4a\md{32}&2\\
     \cline{2-5}   &&a\equiv 1\md4&b\equiv -2^{2k}a\md{2^{2k+3}}&1\\
           \cline{3-5}   & \nu_2(b)=2k (k\ge2)&a\equiv 7\md8&b\equiv 2^{2k}(2-a)\md{2^{2k+2}}&1\\
      \cline{3-5}   & &a\equiv 3\md8&b\equiv 2^{2k}(4-a)\md{2^{2k+3}}&1\\
      \cline{3-5}   & &a\equiv 7\md8&b\equiv -2^{2k}a\md{2^{2k+3}}&2\\    
	     \hline
     a\equiv 0\md{2}&a\equiv 0\md{4}& b\equiv 15-a\md{16}&&1\\
      \cline{2-5} b\equiv 1\md{2}&a\equiv 6\md{8}&\mbox{ See Table A}&&\\
      \hline
       a\equiv 1\md{2}&a\equiv 1\md{8}& b\equiv 1\md{8}&&1\\
      \cline{2-5} b\equiv 1\md{2}&a\equiv 5\md{8}& b\equiv 1\md{8}&&1\\
      \hline
	    
	     {\mbox Otherwise }&&&&0\\
 
   \hline
   \end{array}$$

      $$Table A: a\equiv 6\md8 { \mbox{ and }} b\equiv 1\md2$$  
       $$\begin{array}{|c|c|c|}
	\hline
 {\mbox conditions}& \nu_2(i(K))\\
 \hline
	     a\equiv 6\md{32} { \mbox{ and }}  b\equiv 231-5a \md{256}&1\\
	    \hline
	 \nu_2(a- 6)=2+k,  \nu_2(b-3a+9)\ge 4+2k \mbox{ and } k\ge 3&1\\
	    \hline
 \nu_2(a- 6)=2+k ,  \nu_2(b-3a+9)\ge 4+2j, k\ge 3,   j\le k-2&1\\
 \mbox{ and } (b-3a+9)_2\equiv 3\md4&\\
	    \hline
      \nu_2(a- 22)=2+k,  \nu_2(b-11a+121)\ge 4+2k \mbox{ and } k\ge 3&1\\
	    \hline
 \nu_2(a- 22)=2+k ,  \nu_2(b-11a+121)\ge 4+2j, k\ge 3,   j\le k-2&1\\
 \mbox{ and } (b-11a+121)_2\equiv 3\md4&\\
	    \hline
      a\equiv 22\md{32} { \mbox{ and }}  b\equiv 247+3a \md{256}&1\\
	    \hline
	     a\equiv 14\md{16}{ \mbox{ and }} b\equiv 55+3a \md{64}&1\\
	    \hline
	     {\mbox Otherwise }&0\\
 
   \hline
	     	\end{array}$$
		 		 In particular, when  $\nu_2(i(K))\ge1$, then  $K$ is not monogenic.
		 		 \end{thm}

\smallskip

\begin{thm} \label{npib3}
The prime  integer $3$ divides $i(K)$ if and only if one of the following conditions holds:
\begin{enumerate}
\item
	  $a\equiv -1 \md 3$, $\nu_3(b)=2k $ for some positive integer $k$, and $b_3\equiv 1 \md 3$.
				\item
				$a\equiv b\equiv 1 \md 3$, $\nu_3(a^2-4b)=2k $ for some positive integer $k$ and  $(a^2-4b)_3\equiv 1\md3$.
								\end{enumerate}
				 In particular, if one of these conditions   holds, then $K$ is not monogenic.\\
				In all the above cases we have  $\nu_3(i(K))=1$.
 \end{thm}

\smallskip

\begin{rem}
\begin{enumerate}
 \item 
The field $K$ can be non-monogenic even if the index $i(K)=1$. It suffices to consider the number field $K$ generated by a root of the polynomial  $F(x)=x^4-13$, which is irreducible over $\Q$, as it is $13$-Eisenstein. Since $13=5+8k$ for $k=1$, we conclude  by \cite[section 4.3, p 138]{GR17} that $K$ is  not monogenic. But  $(a,b)=(0,-13)$   satisfies neither the conditions of Theorem \ref{npib2}  nor the conditions of Theorem \ref{npib3},
therefore $i(K)=1$.
\item
The unique method which allows to test whether $ K $ is monogenic is to calculate the solutions of the index form equation of the field $K$ (cf. \cite{G19}, see also Section \ref{444.222}).
\item
It is well known that the index of a quartic field satisfies
$i(K) \in  \{1, 2, 3, 4, 6, 12\}$ (see \cite[p. 234]{En2}).  
Thus for every prime integer $p \ge 5$, $\nu_p(i(K))=0$.
\end{enumerate}
\end{rem}

\vspace{0.5cm}

\subsection{The number of inequivalent generators of power integral bases}

Number fields usually only have a few inequivalent generators of power integral bases
(see e.g. the tables of \cite{G19}). Considering recently some types of monogenic binomials
$x^n-m$ we had the experience that up to equivalence the root of the binomial is 
the only generator of power integral bases in the number field generated by the root.
Calculating "small" solutions (with coefficients $<10^{100}$ in absolute values in the integral bases)
we found that this phenomenon occurs in pure sextic fields generated by a root of
a monogenic binomial $x^6-m$ and in pure octic fields 
generated by a root of 
a monogenic binomial $x^8-m$ both for with $0<m<-5000$.

Therefore we found that it will also 
be interesting to consider the number of inequivalent generators of power 
integral bases of monogenic trinomials of type $x^4+ax^2+b$, utilizing the special properties
of this trinomial.
Note that using \cite{gppsys}, \cite{gppsim} 
in any specific quartic field it is possible to calculate
all generators of power integral bases, but we would like to formulate
here a statement in a parametric form covering an infinite number of fields.
We succeeded to cover the cases $a>1,b>1$:

\begin{thm}
\label{th1}
Assume $a>1,b>1$ and $F(x)=x^4+ax^2+b\in\Z[x]$ is irreducible and monogenic.
If $a,b$ are not of type
\begin{equation}
a=\frac{u\pm 1}{v},\;\;\; b=\frac{u^2-1}{4v^2}
\label{ab}
\end{equation}
for some $u,v\in\Z,v\ne 0,u\ne\pm 1$, then up to equivalence, the root $\al$
of $F(x)$ is the only generator of power integral bases in $K=\Q(\al)$.
\end{thm}

Note that if $a,b$ are of type (\ref{ab}) then there may be several 
inequivalent generators of power integral bases in the corresponding quartic fields. { Hence if $a,b$ are of type (\ref{ab}), then the statement of our Theorem \ref{th1} is not valid.}

We list a few examples of pairs of positive parameters $(a,b)$, represented
in the form (\ref{ab}), such that the trinomial $F(x)=x^4+ax^2+b$ is irreducible and monogenic
but the number field $K=\Q(\alpha)$, generated by a root $\alpha$ of 
$F(x)$ admits several inequivalent generators of power integral bases.
As in such fields $(1,\alpha,\alpha^2,\alpha^3)$ is an integral basis, we shall
give the triples $(x,y,z)$ such that $\gamma=x\alpha+y\alpha^2+z\alpha^3$ generates
a power integral basis in $K$:\\

\noindent
{\bf Case 1.}


{
 $a=\frac{u- 1}{v},\;\; b=\frac{u^2-1}{4v^2}$.

\noindent
$(u,v,a,b)= (3, 1, 2, 2),\;\; (x,y,z)=
                            (1, -1, 0),
                            (1, 0, 0),
                            (1, 0, 1),
                            (1, 1, 0).$\\
$(u,v,a,b)=( 5, 1, 4, 6),\;\; (x,y,z)=
														(1, -1, 0),
                            (1, 0, 0),
                            (1, 1, 0).$\\
$(u,v,a,b)=(  7, 2, 3, 3),\;\; (x,y,z)=
                           ( 1, 0, 0),
                           ( 1, 0, 1),
                           ( 2, -1, 1),
                           ( 2, 1, 1).$\\
$(u,v,a,b)=(  9, 2, 4, 5),\;\; (x,y,z)=
                           ( 1, 0, 0),
                           ( 2, 0, 1).$\\

\noindent
{\bf Case 2.} $a=\frac{u+ 1}{v},\;\; b=\frac{u^2-1}{4v^2}$.

\noindent
$(u,v,a,b)= ( 3, 1, 4, 2),\;\; (x,y,z)=
                           ( 1, 0, 0),
                           ( 3, 0, 1).$\\

\vspace{0.5cm}

\section{A short introduction to Newton polygons}
\label{intro}

We use Newton polygon techniques. This is a standard method which is rather technical but
very efficient to apply. We have introduced the corresponding concepts in several former papers.
Here we only give a brief introduction which makes our proofs understandable.
For a detailed description we refer to Section 3 of our paper \cite{FG8},
to Guardia and Nart \cite{GN} and to Guardia, Montes and Nart \cite{GMN}.

\vspace{0.5 cm}

Let $F(x)\in\Z[x]$ be a monic irreducible polynomial with a root $\al$ and set $K=\Q(\al)$.
We shall use Dedekind's theorem 
\cite[ Chapter I, Proposition 8.3]{Neu}
relating the prime ideal factorization $p\Z_K$ 
and the factorization of $F(x)$ modulo $p$ (for primes $p$ not dividing $(\Z_K:\Z[\al])$).
Also, we shall need Dedekind's criterion \cite[Theorem 6.1.4]{Co} on the divisibility of
$(\Z_K:\Z[\al])$ by primes $p$.

For any  prime integer  $p$, let $\nu_p$ be the $p$-adic valuation of $\Q$, $\Q_p$ its $p$-adic completion, and $\Z_p$ the ring of $p$-adic integers. Let  $\nu_p$  { be } the Gauss's extension of $\nu_p$ to $\Q_p(x)$. 
For any polynomial $P=\sum_{i=0}^na_ix^i\in\Q_p[x]$ we set 
$\nu_p(P)=\mbox{min}(\nu_p(a_i), \, i=0,\dots,n)$.
For nonzero polynomials $P,Q\in\Q_p[x]$, we extend this valuation to 
$\nu_p(P/Q)=\nu_p(P)-\nu_p(Q)$.

Let $\phi\in\z_p[x]$ be a   monic polynomial  whose reduction is irreducible  in
$\F_p[x]$, let $\fph$ be 
the field $\frac{\F_p[x]}{(\overline{\phi})}$. For any
monic polynomial  $F(x)\in \z_p[x]$, upon  the Euclidean division
 by successive powers of $\ph$, we  expand $F(x)$ as 
$F(x)=\sum_{i=0}^la_i(x)\phi(x)^{i},$ called    the $\phi$-{\it expansion} of $F(x)$
 (for every $i$, deg$(a_i(x))<$
deg$(\phi)$). 
The $\ph$-{\it Newton polygon} $\nph{F}$ of $F(x)$ with respect to $p$ is the lower boundary convex envelope of the set of points $\{(i,\nu_p(a_i(x))),\, a_i(x)\neq 0\}$. It is the process of joining the obtained edges  $S_1,\dots,S_r$ ordered by   increasing slopes, which  can be expressed as $\nph{F}=S_1+\dots + S_r$. 
For every side $S_i$ of $\nph{F}$, the length $l(S_i)$ of $S_i$ is the  length of its projection to the $x$-axis and its height  $h(S_i)$  is the  length of its projection to the $y$-axis. We call $d(S_i)=\gcd(l(S_i), h(S_i))$ the  degree of $S$.
 The {\it principal} $\ph$-{\it Newton polygon} $\npp{F}$ of ${F}$
 is the part of the  polygon $\nph{F}$, which is  determined by joining all sides of negative  slopes.
  To every side $S$ {of} $\npp{F}$, with initial point $(s, u_s)$ and length $l$, and to every 
$0\le i\le l$, we attach the following
 residue coefficient $c_i\in\fph$:
$$c_{i}=
\displaystyle\left
\{\begin{array}{ll} 0,& \mbox{ if } (s+i,{\it u_{s+i}}) \mbox{ lies strictly
above } S,\\
\left(\dfrac{a_{s+i}(x)}{p^{{\it u_{s+i}}}}\right)
\,\,
\md{(p,\phi(x))},&\mbox{ if }(s+i,{\it u_{s+i}}) \mbox{ lies on }S,
\end{array}
\right.$$
where $(p,\phi(x))$ is the maximal ideal of $\z_p[x]$ generated by $p$ and $\ph(x)$. 
Let $-\la=-h/e$ be the slope of $S$, where  $h$ and $e$ are two positive coprime integers. Then  $d=l/e$ is the degree of $S$.  Notice that, 
the points  with integer coordinates lying on $S$ are exactly $(s,u_s),(s+e,u_{s}-h),\cdots, (s+de,u_{s}-dh)$. Thus, if $i$ is not a multiple of $e$, then 
$(s+i, u_{s+i})$ does not lie in $S$, and so $c_i=0$. The polynomial
$R_\la(F)(y)=t_dy^d+t_{d-1}y^{d-1}+\cdots+t_{1}y+t_{0}\in\fph[y],$  is  called  
the {\it residual polynomial} of $F(x)$ associated to the side $S$, 
where for every $i=0,\dots,d$,  $t_i=c_{ie}$.

    Let $\npp{F}=S_1+\dots + S_r$ be the principal $\ph$-Newton polygon of $f$ with respect to $p$.     We say that $F$ is a $\ph$-{\it regular polynomial} with respect to $p$, if  $F_{S_i}(y)$ is square free in $\fph[y]$ for every  $i=1,\dots,r$.      The polynomial $F(x)$ is said to be  {\it $p$-regular}  
if $\overline{F(x)}=\prod_{i=1}^r\overline{\ph_i}^{l_i}$ for some monic polynomials $\ph_1,\dots,\ph_t$ of $\Z[x]$ such that $\ol{\ph_1},\dots,\ol{\ph_t}$ are irreducible coprime polynomials over $\F_p$ and    $F(x)$ is  a $\ph_i$-regular polynomial with respect to $p$ for every $i=1,\dots,t$.

Let $\ph\in\Z_p[x]$ be a monic polynomial, with $\overline{\ph(x)}$ irreducible in $\F_p[x]$. As defined in \cite[Def. 1.3]{EMN},   the $\ph$-{\it index} of $F(x)$, denoted by $\ind_{\ph}(F)$, is  deg$(\ph)$ times the number of points with natural integer coordinates that lie below or on the polygon $\npp{F}$, strictly above the horizontal axis,{ and strictly beyond the vertical axis}.

Let $\overline{F(x)}=\prod_{i=1}^r\overline{\ph_i}^{l_i}$ be the factorization of $\overline{F(x)}$ in $\F_p[x]$, where every $\ph_i\in\Z[x]$ is monic polynomial, with $\overline{\ph_i(x)}$ is irreducible in $\F_p[x]$, $\overline{\ph_i(x)}$ and $\overline{\ph_j(x)}$ are coprime when $i\neq j$ and $i, j=1,\dots,t$.
For every $i=1,\dots,t$, let  $N_{\ph_i}^+(F)=S_{i1}+\dots+S_{ir_i}$ be the principal  $\ph_i$-Newton polygon of $F(x)$ with respect to $p$. For every $j=1,\dots, r_i$,  let $F_{S_{ij}}(y)=\prod_{k=1}^{s_{ij}}\psi_{ijk}^{a_{ijk}}(y)$ be the factorization of $F_{S_{ij}}(y)$ in $\F_{\ph_i}[y]$. 
We shall use the following  index theorem of Ore (see e.g. \cite[Theorem 1.7 and Theorem 1.9]{EMN}): 
	
\vspace{0.5cm}
	
 \begin{thm}\label{ore} \mbox{} \\
 \begin{enumerate}
 \item
Using the above notation we have
 $$\nu_p(\ind(F))\ge \sum_{i=1}^r \ind_{\ph_i}(F).$$  The equality holds if $F(x)$ is $p$-regular. 
\item
If  $F(x)$ is $p$-regular, then
$$p\Z_K=\prod_{i=1}^r\prod_{j=1}^{r_i}
\prod_{k=1}^{s_{ij}}\p^{e_{ij}}_{ijk},$$ 
is the factorization of $p\Z_K$ into powers of prime ideals of $\Z_K$ lying above $p$, where $e_{ij}=l_{ij}/d_{ij}$, $l_{ij}$ is the length of $S_{ij}$,  $d_{ij}$ is the degree
 of   $S_{ij}$, and $f_{ijk}=\mbox{deg}(\ph_i)\times \mbox{deg}(\psi_{ijk})$ is the residue degree of the prime ideal  $\p_{ijk}$ over $p$.
 \end{enumerate}
\end{thm}
\begin{rem}\label{remore}
By Dedekind's criterion, $\nu_p(\ind(F))=0$ if and only if $\ind_{\ph_i}(F)=0$ for every $i=1,\dots,r$, which menas $N_{\ph_i}^-(F)$ has a single side of height $1$ for every $i=1,\dots,r$.
\end{rem}
{
When  the program of Ore fails, that is $F(x)$ is not $p$-regular, then in order to complete the factorization of $F(x)$, Guardia, Montes, and Nart introduced the notion of {\it high order Newton polygon}. They showed,  thanks to a theorem  on the index
\cite[Theorem 4.18]{GMN}, that  after a finite number of iterations this process yields all monic irreducible factors of $F(x)$, all prime ideals of $\Z_K$ lying above a prime $p$, the index $(\Z_K:\Z[\al])$, and the absolute discriminant of $K$.  For more details, we refer to \cite{GMN}}.

\vspace{0.5 cm}

\section{Examples}

In this section we give several example that illustrate our main results.

\subsection{Examples 1. Some monogenic number fields defined by non-monogenic trinomials} \label{monon}

In the following statement we give an infinite family of monogenic  number fields generated by roots
of non-monogenic trinomials:

\begin{prop}\label{ppprop}
Let  $K$ be the number field generated by a root $\alpha$ of {an irreducible trinomial } 
\[		
F(x) = x^{4} + 4 a x^2 + 8b  \in \Z[x],
\]
assuming $2$ does not divide $b$, $b$ square-free, and for every odd prime $p$,               
if $p$  does not divide $b$, then $p^2$ does not divide $a^2-2b$ or  $-2^{-1}a\not\in \F_p^2$. \\
Then $F(x)$ is a non-monogenic polynomial, but $K=\Q(\al)$ is a monogenic number field  {and $\theta=\frac{\al^3}{4}$ is a generator of a power integral basis of $K$.}
\end{prop}

For the proof of this Proposition see Section \ref{sectionproof}.

\vspace{0.5cm}

\subsection{Examples 2: Monogenic and non-monogenic fields, index forms.}	
\label{444.222}

	Now let $K$ be a number field generated by a root $\alpha$ of $x^4+ax^2+b\in\Z[x]$ such that { $b$ is square-free $a\equiv 0\md2$, $b\equiv 1\md2$, and for every odd prime $p$ either $p^2$ does not divide $a^2-4b$ 
	or $-2^{-1}a\not\in \F_p^2$.} Then

\smallskip
\noindent
(1)
{	If $\nu_2(b+1+a)=1$, then by Theorem \ref{intclos}, $\Z[\al]$ is the ring of integers of $K$.}\\

\smallskip
\noindent
(2)
	 {If $\nu_2(b+1+a)\ge 3$ and $a\equiv 0\md4$, then by Theorem \ref{npib2} (9), $2$ divides $i(K)$, and so $K$ is not mongenic.}\\

\smallskip
\noindent
(3)
	If $\nu_2(b+1+a)=2$  and $a\equiv 0\md4$, then $\ol{F(x)}=\ph^4$, where $\ph=x-1$.
	The principal $\ph$-Newton polygon (cf. \cite{FG8})
	of $F$ with respect to $2$ is  
	$\nph{F}=S$ having 
	a single side joining $(0,2)$, $(2,1)$, and $(4,0)$. Thus $S$ is of degree $2$, slope $\la=-1/2$, and $R_{1/2}(F)(y)=y^2+y+1$ is irreducible over $\fph$. Thus $\left(1,\al,\frac{\al^2+1}{2}, \frac{\al^3+\al}{2}\right)$ is a $\Z$-basis of $\Z_K$.

The conditions $\nu_2(b+1+a)=2$  and $a\equiv 0\md4$ imply $a=4k,b=4l-4k-1$ with integer parameters $k,l$.
The index form corresponding to the above integral basis is $F_1\cdot F_2$ where

{\small
\[
F_1(x_1,x_2,x_3)=-x_1^2-x_3^2l+2x_3^2k-x_1x_3+2x_1x_3k,	
\]
\[
F_2(x_1,x_2,x_3)=
4x_2^4k^2-4x_3^4l-4x_2^4l+4x_2^4k+64x_3^4k^4-32x_3^4k^3+4x_3^4l^2+8x_1^3x_3+8x_1^2x_3^2+4x_1x_3^3-2x_3^2x_2^2
\]
\[
-4x_3^4k
+20x_3^4k^2+x_2^4+x_3^4+4x_1^4-32x_1^3x_3k-128x_1x_3^3k^3+8x_1^2x_2^2k-24x_1x_3^3k+8x_2^2x_3^2k^2
\]
\[
+64x_1x_3^3k^2-40x_1^2x_3^2k
-8x_1x_3^3l+32x_2^2x_3^2k^3+8x_3^4kl+8x_3^2x_2^2l-32x_3^4k^2l-8x_1^2x_3^2l+96x_1^2k^2x_3^2
\]
\[
-4x_1x_2^2x_3+32x_1x_3^3kl-24x_2^2x_3^2kl-8x_1x_2^2x_3k+16x_3x_2^2x_1l-32x_1x_2^2x_3k^2.
\]
}
These number fields can be either monogenic or not.\\

\smallskip
\noindent
(3.1) For $k=l=1$ we obtain the number field $K=\Q(\alpha)$ generated by a root $\alpha$ 
of $f(x)=x^4+4x^2-1$. The field $K$ is a mixed quartic field with Galois group $D_4$
and discriminant $d_K=-400$. 
Up to equivalence all generators of power integral bases are given by the solutions
 $(x_1,x_2,x_3)=(1,0,1),(2,0,1),(3,-1,2),(3,1,2)$ of the index form equation.
(Here and in the following we used the algorithm described in Ga\'al, Peth\H o and Pohst 
\cite{gppsim} (see also \cite{G19}) to solve the index form equations in order to find all 
generators of power integral bases of the corresponding quartic fields.)
\\

\smallskip
\noindent
(3.2) For $k=2,l=7$ we obtain the number field $K=\Q(\alpha)$ generated by a root $\alpha$ 
of $f(x)=x^4+8x^2+19$. The field $K$ is a totally complex quartic field with Galois group $D_4$
and discriminant $d_K=2736$, which is not monogenic, since the index form equation
has no solutions.\\

\smallskip
\noindent
(4)
	{If $\nu_2(b+1+a)=3$  and $a\equiv 2\md4$, then $\ol{F(x)}=\ph^4$, where $\ph=x-1$, 
		The principal $\ph$-Newton polygon	of $f$ with respect to $2$ is	
	$\nph{F}=S$ has a single side joining $(0,3)$ and $(4,0)$. Thus $S$ is of degree $1$, and so 
	the residual polynomial $R_{\la}(F)(y)$ (cf. \cite{FG8}) is irreducible over $\fph$. Thus 
	$\left(1,\al,\frac{\al^2+1}{2}, \frac{\al^3+\al^2+3\al+3}{4}\right)$ is a $\Z$-basis of $\Z_K$.} 
	
The conditions $\nu_2(b+1+a)=3$  and $a\equiv 2\md4$ imply $a=4k+2,b=8l-4k-3$ with integer parameters $k,l$.
The index form corresponding to the above integral basis is $F_1\cdot F_2$ where

{\small
\[
F_1(x_1,x_2,x_3)=-2x_1^2-x_3^2l+2x_3^2k-2x_1x_3+2x_1x_3k,
\]
\[
F_2(x_1,x_2,x_3)=4x_1^4+4x_1^3x_3+4x_1^2x_3^2+16x_3x_2^2x_1l-12x_3x_2^2x_1k+16x_3^2x_2x_1l-12x_3^2x_2x_1k+8x_1x_3^3kl
\]
\[
-12x_3^2x_2^2kl-12x_3^3x_2kl-16x_3x_2^2x_1k^2-16x_3^2x_2x_1k^2+8x_3x_2x_1^2k+4x_1x_3^3k^2+8x_3x_2^3k^2
\]
\[
+24x_1^2x_3^2k^2-2x_3^4kl-4x_1^2x_3^2l-4x_3^4k^2l-16x_1^3x_3k+4x_3^3x_2k+8x_3^2x_2^2k^3+2x_1x_3^3l
\]
\[
-8x_1^2x_3^2k+12x_3^2x_2^2k+2x_3^3x_2l-6x_3^2x_2^2l+16x_3x_2^3k+8x_3^3x_2k^2+8x_3^3x_2k^3-8x_1x_3^3k
\]
\[
+8x_2^2x_1^2k-16x_1x_3^3k^3-16x_3x_2^3l+12x_3^2x_2^2k^2+4x_3^4k^2+4x_3^4k^4+4x_2^4k^2+8x_2^4k+x_3^4l^2
\]
\[
-8x_2^4l+4x_2^4+4x_3^2x_2^2+8x_3x_2^3+4x_2^2x_1^2+4x_3x_2x_1^2-4x_3x_2^2x_1-4x_3^2x_2x_1,	
\]
}
These number fields can be either monogenic or not.\\

\smallskip
\noindent
(4.1) For $k=l=1$ we obtain the number field $K=\Q(\alpha)$ generated by a root $\alpha$ 
of $f(x)=x^4+6x^2+1$. The field $K$ is a totally complex quartic field with Galois group $V_4$
and discriminant $d_K=256$. 
Up to equivalence all generators of power integral bases are given by the solutions
 $(x_1,x_2,x_3)=(1,-1,1),(1,0,1)$ of the index form equation.\\

\smallskip
\noindent
(4.2) For $k=1,l=3$ we obtain the number field $K=\Q(\alpha)$ generated by a root $\alpha$ 
of $f(x)=x^4+6x^2+17$. The field $K$ is a totally complex quartic field with Galois group $D_4$
and discriminant $d_K=4352$, which is not monogenic, since the index form
equation has no solutions.\\

\vspace{0.5cm}

\subsection{Examples 3: Applying Engstrom's results.}

Let $K$ be a number field generated by a root of the irreducible polynomial $F(x)=x^4+ax^2+b\in \Z[x]$. In the following examples we show that $i(K)\neq 1$, which implies that $K$ is not monogenic. 
In all these examples we shall use the result of Engstrom \cite[p. 234]{En2} stating that the exponents 
of 2 and 3 in the field index $i(K)$ of $K$ only depend on the type of factorizations
of $2\Z_K$ and $3\Z_K$ into prime ideals of $\Z_K$.

	\begin{enumerate}
		\item 
				If $a=-55$ and $b=-99$, then  $F(x)$ is $11$-Eisenstein, and so $F(x)$ is irreducible over $\Q$. Since $a\equiv 1\md8$ and $b\equiv 5\md8$, by Theorem \ref{npib2}(1) and its proof, $2$ divides $i(K)$ and $2\Z_K=
				\p_1\p_2$ with residue degree $2$ each prime factor. By Engstrom's theorem, we conclude that $\nu_2(i(K))=1$. Again by Theorem \ref{npib3}(1) and its proof, $3$ divides $i(K)$ and
				$3\Z_K=
				\p_1\p_2\p_3\p_4$ with residue degree $1$ each prime factor. By Engstrom's theorem, we conclude that 
				$\nu_3(i(K))=1$. Hence  $i(K)=6$.		
				\item
					Similarly, if  $a=17$ and $b=765$, then  $F(x)$ is  irreducible over $\Q$. Since  $2\Z_K=
				\p_1\p_2$ with residue degree $2$ each prime factor and 
				$3\Z_K=
				\p_1\p_2\p_3\p_4$ with residue degree $1$ each prime factor,  we conclude by Engstrom's theorem, that $\nu_2(i(K))=1$, $\nu_3(i(K))=1$. Hence  $i(K)=6$.		 
						\item
					If  $a=3$ and $b=4$,  then  $F(x)$ is irreducible over $\Q$. Since $\nu_2(b)=2$ and $a+b_2\equiv 4\md8$, we conclude that $2\Z_K=\p_1\p_2\p_3\p_4$ with residue degree $1$ each prime factor. By Engstrom's theorem it implies $\nu_2(i(K))=2$, For $p=3$, since $3$ divides $a$ and does not divide $b$, hence $3$ does not divide $(\Z_K:\Z[\al])$. 				
					Thus by $\nu_3(i(K))=0$, $\nu_2(i(K))=2$ we have $i(K)=4$.
					\item
					If  $a=-1$ and $b=304$,  then  $F(x)$ is irreducible over $\Q$. Since $\nu_2(b)=4$ and $1+a+b\equiv 0\md8$, we conclude that $2\Z_K=\p_1\p_2\p_3\p_4$ with residue degree $1$ each prime factor.   For $p=3$, since $3$ divides $a\equiv -1\md3$, $\nu_3(b)=2$, and $b_3\equiv 1\md3$, we have $3\Z_K=\p_1\p_2\p_3\p_4$ with residue degree $1$ each prime factor. By Engstrom's theorem we have $\nu_3(i(K))=1$, $\nu_2(i(K))=2$, hence  $i(K)=12$.
				
					\item
					If  $a=4$ and $b=-86$,  then  $F(x)$ is $2$-Eisenstein, and so $F(x)$ is irreducible over $\Q$ and $2$ does not divide $(\Z_K:\Z[\al])$. Since $a\equiv 1\md3$ and $b+1+a\equiv 0\md{27}$,  we conclude that $3\Z_K=
				\p_1\p_2\p_3\p_4$ with residue degree $1$ each prime factor. By Engstrom's theorem it implies $\nu_3(i(K))=1$. By $\nu_2(i(K))=0$, $\nu_3(i(K))=1$
				we have $i(K)=3$.
					\item
					If  $a=-1$ and $b=304$,  then  $F(x)$ is irreducible over $\Q$. Since $\nu_2(b)=4$ and $1+a+b\equiv 0\md8$, we conclude that $2\Z_K=\p_1\p_2\p_3\p_4$ with residue degree $1$ each prime factor.  For $p=3$, since $3$ divides $a\equiv -1\md3$, $\nu_3(b)=2$, and $b_3\equiv 1\md3$, we have $3\Z_K=\p_1\p_2\p_3\p_4$ with residue degree $1$ each prime factor. By Engstrom's theorem we have  $\nu_3(i(K))=1$, $\nu_2(i(K))=2$, hence  $i(K)=12$.
					\end{enumerate}

\section{Proofs of our main results}
\label{sectionproof}

\noindent
{\bf Proof of {Theorem}
 \ref{intclos}}.
\\
Let $p$ be a prime integer candidate to divide  $(\Z_K:\Z[\al])$. 
Since $\triangle=16b(a^2-4b)^2$ is the discriminant of $F(x)$,  we conclude that  $p=2$, or $p$ divides $b$ or $p$ divides $a^2-4b$.  
\begin{enumerate}
\item
If $p$ divides both  of $a$ and $b$, then $\overline{F(x)}=\ph^4$, where $\ph=x$. If $\nu_p(b)=1$, then by Remark \ref{remore},
 $\nu_p((\Z_K:\Z[\al]))=0$ if and only if  $\nu_p(b)=1$.
\item
For $p=2$, 
\begin{enumerate}
\item
If $2$ divides $a$ and $b$, then by the first point, $2$ does not divide  $(\Z_K:\Z[\al])$ if and only if $\nu_2(b)=1$.
 \item
 If $2$ divides $b$ and does not $a$, then  $\ol{F(x)}=x^2(x+1)^2$. Let $\ph=x-1$. Then
 $F(x)=\ph^4+4\ph^3+(6+a)\ph^2+(4+2a)\ph+(1+b+a)$. Thus, by Remark \ref{remore},  we conclude that  $2$ does not divide  $(\Z_K:\Z[\al])$ if and only if ind$_x(F)=0$ and ind$_\ph(F)=0$, which means that $\nu_2(b)=1$ and $\nu_2(b+a+1)=1$. That is  $b\equiv 2\md4$ and $a\equiv 2-(1+b)\equiv 3 \md4$.
 \item
 If $2$ divides $a$ and does not divide $b$, then  $\ol{F(x)}=(x+1)^4$. Let $\ph=x-1$. Then
 $F(x)=\ph^4+4\ph^3+(6+a)\ph^2+(4+2a)\ph+(1+b+a)$. By  Remark \ref{remore}, we conclude that  $2$ does not divide  $(\Z_K:\Z[\al])$ if and only if  $\nu_2(b+a+1)=1$, which means  $b\equiv 1-a \md4$.
 \item
 If $2$ does not divide $ab$, then  $\ol{F(x)}=(x^2+x+1)^2$. Let $\ph=x^2+x+1$ and $F(x)=\ph^2+(a-1-2x)\ph+(1-a)x+b-a$. Thus by Remark \ref{remore},  $2$ does not divide  $(\Z_K:\Z[\al])$ if and only if $\nu_2(a-1)=1$ or $\nu_2(a-b)=1$, which means $a\equiv 3\md4$ or $b\equiv 2+ a \md4$.
 \end{enumerate}
 \item
If $p$ is odd, $p$ divides $b$ and does not divide $a$, then   $\ol{F(x)}=x^2(x^2+a)$. Since $x^2+a$ is square free in $\F_p[x]$, then by Remark \ref{remore},   $p$ does not divide  $(\Z_K:\Z[\al])$ if and only if ind$_x(F)=0$. That is $\nu_p(b)=1$.
\item
{Now assume that $p$ is odd and $p$ does not divide $ab$. If $p$ does not divide $a^2-4b$, then  $p^2$ does not divide $\triangle(F)$, and so $p$ does not divide $(\Z_K:\Z[\al])$. Also since $F'(x)=2x(2x^2+a)$, if  $-2^{-1}a\not\in \F_p^2$, then  $\ol{F(x)}$ is square free in $\F_p[x]$. Thus, by Dedekind's criterion, we conclude that  $p$ does not divide  \mbox{$(\Z_K:\Z[\al])$}. If  $p$  divides $a^2-4b$ and $-2^{-1}a\in \F_p^2$, then by Hensel's lemma, let $u\in \Z_p$ such that $2u^2+a\equiv 0\md{p^r}$ with $r\ge 2$ a positive integer. Then $2^2F(u)=(-a)^2+2a\cdot (-a)+4b\equiv -(a^2-4b)\md{p^r}$. Let $\ph=x-u$ and $F(x)=\ph^4+a_3\ph^3+a_2\ph^2+a_1\ph+a_0$ be the $\ph$-expansion of $F(x)$ with $a_1= F'(u)$ and $a_0= F(u)$. Thus $a_1\equiv F'(u)\equiv 0\md{p^2}$ and $a_0\equiv 0\md{p}$. By Remark \ref{remore}, $ind(F)=0$ if and only if $\nu_p(F(u))=1$, which means  that $p^2$ does not divide $a^2-4b$.}
\end{enumerate}
$\hfill\Box$

		\smallskip
	
	For the proof of Theorems \ref{npib2} and \ref{npib3}, based on Engstrom's results, we need to factorize $2\Z_K$ and $3\Z_K$ into the product of powers of prime ideals of $\Z_K$. Recall the following lemma, which  characterizes the prime  index divisors of $K$. Its proof is an immediate consequence of Dedekind's theorem:

	 	\begin{lem}\label{index}
		Let $p$ be a rational prime integer and $K$ be a number field. For every positive integer $f$, let $\mathcal{P}_f$ be the number of distinct prime ideals of $\mathbb{Z}_K$ lying above $p$ with residue degree $f$ and $\mathcal{N}_f$ the number of monic irreducible polynomials of $\mathbb{F}_p[x]$ of degree $f$.  Then $p$ is a prime common index divisor of $K$ if and only if $\mathcal{P}_f>\mathcal{N}_f$ for some positive integer $f$.
	\end{lem}

\smallskip

\noindent
{\bf Proof of Theorem \ref{npib2}}.
 {By virtue of Engstrom's results \cite{En2}, we need to factorize $2\Z_K$ into powers of prime ideals of $\Z_K$. Moreover, $2$ divides $i(K)$ if and only if  $2\Z_K=\p_1\p_2\p_3\p_4$ with residue degree $1$ each prime ideal  or $2\Z_K=\p_1\p_2$ with residue degree $2$ each prime ideal or  also $2\Z_K=\p_1^2\p_2\p_3$ with residue degree $1$ each prime ideal.}  { According to Theorem \ref{intclos}, we deal only with the   cases, where $\Z[\al]$ is not integrally closed.}
\begin{enumerate}
	 \item 
	  If $\nu_2(b)\ge 2$ and $\nu_2(a)\ge 1$,   then  $ \ol{F(x)}= \ph^4$ in $\F_2[x]$, where $\ph = x$.
	If $\nph{F}=S$ has a single side, then let $d$ be the degree of $S$. Since $\nu_2(a)\le 1$ or $\nu_2(b)\le 3$, we conclude that $d\in\{1,2\}$. If $d=1$, then there is a unique prime ideal of $\Z_K$ lying above $2$ with residue degree $1$. If  $d=2$ ($\nu_2(a)\ge 1$ and $\nu_2(b)= 2$), then there is a unique prime ideal of $\Z_K$ lying above $2$ with residue degree $2$ and ramification index $2$ or  two prime ideals of $\Z_K$ lying above $2$ with residue degree $1$ and ramification index $2$ each. 	If $\nph{F}=S_1+S_2$ has two sides, that is $\nu_2(a)=1$ and $\nu_2(b)\ge 3$, then there are two cases: If $\nu_2(b)$ is even then  there are two prime ideals of $\Z_K$ lying above $2$ with residue degree $1$ and ramification index $2$ each.  If $\nu_2(b)=2k+1$ for some positive integer $k$, then let $\ph=x+2^k$ and $F(x)=\ph^4+2^{k+2}\ph^3+(a+3\cdot 2^{2k+2})\ph^2+(2^{k+1}a+2^{3k+2})\ph+(2^{4k}+2^{2k}a+b)$ be the $\ph$-expansion of $F(x)$. Since   $\nu_2(a+3\cdot 2^{2k+2})=1$, $\nu_2(2^{k+1}a+2^{3k+2})=k+2$ and   $\nu_2(2^{4k}+2^{2k}a+b)\ge  2k+3$, we conclude that $2$ divides $i(K)$ if and only if $\nph{F}$ has three sides. That is $\nu_2(2^{4k}+2^{2k}a+b)\ge  2k+4$, which means that  $k\ge 2$ and $\nu_2(2^{2k}a+b)\ge 2k+4$ or $k=1$ and $\nu_2(2^4+4a+b)=6$ ($b\equiv 48-4a\md{64}$). In these cases, there are three prime ideals of $\Z_K$ lying above $2$ with residue degree $1$  each.  Thus  $\nu_2(i(K))=1$.
 \item
If  $\nu_2(b)\ge 1$ and   $a\equiv 1\md2$, then  $ \ol{F(x)}= (x\cdot \ph)^2$ in $\F_2[x]$, where  $\ph = x-1$. 
 We have the following cases:
\begin{enumerate}
\item    
  If   $\nu_2(b)=1+2k$ for some positive integer $k$, then $x$  provides  a unique prime ideal of $\Z_K$ lying above $2$ with residue degree $1$. In this case $2$ is a common index divisor of $K$ if and only if  $\ph$ provides two prime ideals of $\Z_K$ lying above $2$ with residue degree $1$ each. For this reason, 
  let   
  $F(x)=\ph^4+4\ph^3+(a+6)\ph^2+(2a+4)\ph+(b+a+1)$. 
    Since $\nu_2(2a+4)=1$, we conclude that $2$ is a common index divisor of $K$ if and only if  $\nu_2(b+a+1)\ge 3$. That is $\nu_2(b)=1+2k$ and  $a\equiv 7-b\md8$. In this case $\nu_2(i(K))=1$.
  \item    
  If   $\nu_2(b)=2k$ for some positive integer $k$, then $N_x^-(F)=S$ has a single side joining $(0,2k)$ and $(2,0)$.   Let $\ph_2=x+2^k$ and  
  $F(x)=\ph_2^4-2^{k+2}\ph_2^3+(a+3\cdot 2^{2k+1})\ph_2^2-
  (2^{k+1}a+2^{3k+2})\ph_2+
  (b+2^{2k}a+2^{4k})$.
    Since $\nu_2(2^{k+1}a+2^{3k+2})=k+1$ and $\nu_2(b+2^{2k}a+2^{4k})\ge 2k+1$, we conclude that
     $2$ is a common index divisor of $K$ if and only if 
  $\nu_2(b+a+1)\ge 3$ and  $\nu_2(b+2^{2k}a)=2k+1$   or  $\nu_2(b+a+1)=1$ and  $\nu_2(b+2^{2k}a)\ge 2k+3$ or
  $\nu_2(b+a+1)=2$ and  $\nu_2(b+2^{2k}a)=2k+2$   or also $\nu_2(b+a+1)\ge 3$ and $\nu_2(b+2^{2k}a))\ge 2k+3$  ({In this last case $\nu_2(i(K))=2$}). \\
   {Remark that if $k=1$, then
      $$\nu_2(b+a+1)\ge 3 \mbox{ and } \nu_2(b+2^{2}a)=3 \Longleftrightarrow a\equiv 3\md8 \mbox{ and } b\equiv 4(2-a)\md{16},$$ 
      $$\nu_2(b+a+1)=1 \mbox{ and } \nu_2(b+2^{2}a)\ge 5 \Longleftrightarrow a\equiv 1\md4 \mbox{ and } b\equiv -4a\md{32},$$ 
    $$\nu_2(b+a+1)=2 \mbox{ and } \nu_2(b+2^{2}a)=4 \Longleftrightarrow a\equiv 7\md8 \mbox{ and }  b\equiv 4(4-a)\md{32}  \mbox{ and }$$
    $$\nu_2(b+a+1)\ge 3 \mbox{ and } \nu_2(b+2^{2}a)\ge 5 \Longleftrightarrow a\equiv 3\md8 \mbox{ and } b\equiv -4a\md{32}.$$
    For $k\ge 2$, then   
   $$\nu_2(b+a+1)\ge 3 \mbox{ and } \nu_2(b+2^{2}a)=2k+1 \Longleftrightarrow a\equiv 7\md8 \mbox{ and } b\equiv 2^{2k}(2-a)\md{2^{2k+2}},$$ 
   $$\nu_2(b+a+1)=1 \mbox{ and } \nu_2(b+2^{2}a)\ge 2k+3 \Longleftrightarrow a\equiv 1\md4 \mbox{ and } b\equiv -2^{2k}a\md{2^{2k+3}},$$ 
    $$\nu_2(b+a+1)=2 \mbox{ and } \nu_2(b+2^{2k}a)=2k+2 \Longleftrightarrow a\equiv 3\md8 \mbox{ and }  b\equiv 2^{2k}(4-a)\md{2^{2k+3}}  \mbox{ and }$$
    $$\nu_2(b+a+1)\ge 3 \mbox{ and } \nu_2(b+2^{2k}a)\ge 2k+3 \Longleftrightarrow a\equiv 7\md8 \mbox{ and } b\equiv -2^{2k}a\md{2^{2k+3}}.$$}
  \end{enumerate}
  \item
  If  $\nu_2(b)=0$ and $\nu_2(a)\ge 1$, then $ \ol{F(x)}=  \ph^4$ in $\F_2[x]$, where  $\ph = x-1$. Let 
  $F(x)=\ph^4+4\ph^3+(a+6)\ph^2+(2a+4)\ph+(b+a+1)$.
   \begin{enumerate}
 \item  
 For $\nu_2(a)\ge 2$,  we have $\nu_2(a+6)=1$ and $\nu_2(2a+4)=2$. It follows that
  $2$ is a common divisor of $K$ if and only if $\nu_2(b+a+1)\ge 4$. That is $b\equiv 15-a \md{16}$. In this case $\nu_2(i(K))=1$.
      \item
For   $a\equiv 2\md8$, let $C=(b+a+1)$, $B=(2a+4)$ and $A=(a+6)$. According to   $a\equiv 2\md{16}$ or $a\equiv 10\md{16}$, we get  $\nu_2(B)=3$ and $\nu_2(A)\ge 3$. It follows that if $\nu_2(a+b+1)\in \{1,3\}$, then $\npp{F}=S$ has a single side of degree $1$. Thus, there is a single prime ideal of $\Z_K$ lying above $2$.
If $\nu_2(a+b+1)=2$, then $\npp{F}=S$ has a single side of degree $2$ and $R_\la(F)(y)=(y+1)^2$. Thus, we have to use second order Newton's polygon techniques. Since $-\la=-1/2$ is the slope of $S$, we conclude that $2$ divides the ramification index $e(\p)$ of $\p$ over $2$ for every  prime ideal $\p$  of $\Z_K$ lying above $2$. Thus the factorization of $2\Z_K$ has one of these forms: $\p_1^4$, $\p_1^2\p_2^2$ or $\p_3^2$ with $f(\p_1)=f(\p_2)=1$ and  $f(\p_3)=2$ are the residue degrees. In all these cases, $2$ is not a common  divisor of $K$.
If $\nu_2(a+b+1)=4$, then $\npp{F}=S$ has a single side of degree $4$ and $R_\la(F)(y)=y^4+y+1$ is irreducible over $\fph=\F_2$. So, there is a single prime ideal of $\Z_K$ lying above $2$. 
If $\nu_2(a+b+1)\ge 5$, then $\npp{F}=S_1+S_2$ has two sides such that $S_1$ is  of degree $1$, $S_2$ is  of degree $3$ and $R_{\la_2}(F)(y)=y^3+1=(y+1)(y^2+y+1)$ in  $\fph[y]$. Hence  there are exactly two prime ideals of  $\Z_K$ lying above $2$ with residue degree $1$ each and a single prime ideal of $\Z_K$ lying above $2$ with residue degree $2$. So, $\nu_2(i(K))=0$. 
 \item
 { For $a\equiv 6\md8$,      
  we have $F(x)=\ph^4+4\ph^3+(a+6)\ph^2+(2a+4)\ph+(b+a+1)$ with $\ph=x-1$. If $\nu_2(b+a+1)\in\{1,3\}$, then $2\Z_K=\p^4$ with $\p$ the prime ideal of $\Z_K$ lying above $2$ with residue degree $1$.  In this case $2$ does  not divide $i(K)$.\\
  If $\nu_2(b+a+1)=2$, then $\nph{F}=S$ has a single side with $R_\la(F)(y)=(y+1)^2$. Thus, there  are two cases, which are: Either there are   
  two prime ideals of $\Z_K$ lying above $2$ with residue degree $1$ and ramification index $2$ each. Or there is   
  a unique prime ideal of $\Z_K$ lying above $2$ with residue degree $2$ and ramification index $2$. Again in this case $2$ does  not divide $i(K)$.
  \item 
 If $\nu_2(b+a+1)\ge 4$, then let $G(x)=x^4+2x^3+Ax^2+Bx+C$ be the minimal polynomial of $\theta=\frac{\al-1}{2}$ over $\Q$, where $A=\frac{a+6}{4}$, $B=\frac{2a+4}{8}$, and $C=\frac{1+b+a}{16}$. Since $a\equiv 6\md8$ and $\nu_2(b+a+1)\ge 4$, $G(x)\in \Z[x]$, $\theta\in \Z_K$ is a primitive element of $K$,  $A\equiv 1\md2$ and   we can replace $F(x)$ by $G(x)$. 
     Since  $a\equiv 6\md8$, then $\nu_2(A)=\nu_2(a+6)-2=0$ and $\nu_2(2a+4)\ge 4$. So, if $\nu_2(b+a+1)=4$, then $\ol{G(x)}=(x^2+x+1)^2$.
    Let $\ph=x^2+x+1$ and $G(x)=\ph^2 + \frac{a-6}{4}\ph+\frac{b-3a+9}{16}$. It follows that:
     \begin{enumerate}
\item 
 If $\nu_2(b-3a+9)=5$, then $2$ does not divide ind$(F)$, and so   $\nu_2(i(K))=0$.  
\item
If $\nu_2(b-3a+9)=6$ and $a\equiv 6\md{32}$, then for $\ph=x^2+x-1$, $G(x)=\ph^2+\frac{a+10}{4}+\frac{b+5a+25}{16}$. Since   $b+5a+25\equiv b-3a+9+8a+16\equiv 8(a+10)\md{128}$, we conclude that  $\nu_2(\frac{b+5a+25}{16})\ge 3$. As $\nu_2(\frac{a+10}{4})=2$, then  $\nu_2(i(K))=1$ if and only if $\nu_2({b+5a+25})\ge 8$, which means $b\equiv 231-5a\md{256}$.
\item
If $\nu_2(b-3a+9)=7$ and $a\equiv 6\md{16}$, then for $\ph=x^2+x+1$, $\npp{G}$ has a single side of degree $1$ and so, $\nu_2(i(K))= 0$.
\item
If $\nu_2(b-3a+9)\ge  8$ and $a\equiv 6\md{32}$, then let $k=\nu_2(a-6)-2$ and $v=\nu_2(b-3a+9)-4$. If $2k\le v$,  then $\npp{G}$ has two sides, and so  $2\Z_K=\p_1\p_2$ with residue degree $2$ each. Therefore,  $\nu_2(i(K))=1$. \\
If $v< 2k$ and $v$ is odd, then $\npp{G}$ has a single side of degree $1$,   $2\Z_K=\p_1^2$, and so  $\nu_2(i(K))=0$. \\
If $v=2k$, then $\npp{G}$ has a single side with  $R_\la(G)(y)=y^2+y+1=(y-x)(y-x^2)$, $2\Z_K=\p_1\p_2$ with residue degree $2$ each, and so   $\nu_2(i(K))=1$. \\
If $v=2j$ for some positive integer $j$ and $v< 2k$, then let $\ph=x^2+x+1+2^j$ and $G(x)=\ph^2+ (\frac{a-6}{4}+2^{j+1})\ph+\frac{b-3a+9-(a-6)2^{j+2}+2^{2j+4}}{16}$.  So,
If  $j<k-1$, then $\nu_2(\frac{a-6}{4}+2^{j+1})=j+1$ and \\
$\nu_2(\frac{b-3a+9-(a-6)2^{j+2}+2^{2j+4}}{16})\ge 2j+1$. Analogously to the previous point, $\nu_2(i(K))=1$ if and only if \\
$\nu_2(\frac{b-3a+9-(a-6)2^{j+2}+2^{2j+4}}{16})\ge 2j+2$. Say $(b-3a+9)_2\equiv 3 \md4$.\\
If  $k=j+1$, then   $\nu_2(\frac{a-6}{4}+2^{j+1})\ge j+1$ and \\
$\nu_2(\frac{b-3a+9-(a-6)2^{j+2}+2^{2j+4}}{16})= 2j-1$ is odd. Thus,  $\npp{G}$ has a single side of degree $1$,   $2\Z_K=\p_1^2$, and so  $\nu_2(i(K))=0$.
\item
If $\nu_2(b-3a+9)=6$ and $a\equiv 22\md{32}$, then for $\ph=x^2+x+3$, $G(x)=\ph^2+\frac{a-22}{4}+\frac{b-11a+121}{16}$.
let $k=\nu_2(a-22)-2$ and $v=\nu_2(b-11a+121)-4$. If $2k\le v$,  then $\npp{G}$ has two sides, and so  $2\Z_K=\p_1\p_2$ with residue degree $2$ each. Therefore,  $\nu_2(i(K))=1$. \\
If $v< 2k$ and $v$ is odd, then $\npp{G}$ has a single side of degree $1$,   $2\Z_K=\p_1^2$, and so  $\nu_2(i(K))=0$. \\
If $v=2k$, then $\npp{G}$ has a single side with  $R_\la(G)(y)=y^2+y+1=(y-x)(y-x^2)$, $2\Z_K=\p_1\p_2$ with residue degree $2$ each, and so   $\nu_2(i(K))=1$. \\
If $v=2j$ for some positive integer $j$ and $v< 2k$, then let $\ph=x^2+x+3+2^j$ and \\
$G(x)=\ph^2+ (\frac{a-22}{4}+2^{j+1})\ph+\frac{b-11a+121-(a-22)2^{j+2}+2^{2j+4}}{16}$.  So,
if  $j<k-1$, then $\nu_2(\frac{a-22}{4}+2^{j+1})=j+1$ and\\
$\nu_2(\frac{b-11a+121-(a-22)2^{j+2}+2^{2j+4}}{16})\ge 2j+1$. Analogously to the previous point, $\nu_2(i(K))=1$ if and only if \\$\nu_2(\frac{b-11a+121-(a-22)2^{j+2}+2^{2j+4}}{16})\ge 2j+2$, say $(b-11a+121)_2\equiv 3 \md4$.\\
If  $k=j+1$, then   $\nu_2(\frac{a-22}{4}+2^{j+1})\ge j+1$ and \\
$\nu_2(\frac{b-3a+9-(a-22)2^{j+2}+2^{2j+4}}{16})= 2j-1$ is odd. Thus,  $\npp{G}$ has a single side of degree $1$,   $2\Z_K=\p_1^2$, and so  $\nu_2(i(K))=0$. 
 \item
If $\nu_2(b-3a+9)\ge  8$ and $a\equiv 22\md{32}$, then for $\ph=x^2+x+1$, $G(x)=\ph^2+\frac{a-6}{4}+\frac{b-3a+9}{16}$.   Since $\nu_2(a-6)-2=2$ and $v=\nu_2(b-3a+9)-4\ge 4$, then  $2\Z_K=\p_1\p_2$ with residue degree $2$ each, and so   $\nu_2(i(K))=1$. \\
\item 
If $a\in\{14, 30\}\md{32}$, then $\nu_2( \frac{a-6}{4})=1$. 
It follows that if $\nu_2(b-3a+9)=6$, then $\npp{G}$ has a single side with  $R_\la(G)(y)=y^2+y+1=(y-x)(y-x^2)$. Therefore $2\Z_K=\p_1\p_2$ with residue degree $2$ each, and so   $\nu_2(i(K))=1$. Finally, if $\nu_2(b-3a+9)\ge 7$,
 then $\npp{G}$ has two sides with degree $1$ each. Thus $2\Z_K=\p_1\p_2$ with residue degree $2$ each, and so   $\nu_2(i(K))=1$.
         \end{enumerate}}
\end{enumerate}
  \item 
	If	$\nu_2(ab)=0$, then  $ \ol{F(x)}= \ph^2$ in $\F_2[x]$, where $\ph = x^2+x+1$. In this case $2$ is a  divisor of $K$ if and only if $\ph$ provides two prime ideals of $\Z_K$ lying above $2$ with residue degree $2$ each.
 Let 	$F(x)= \ph^2+(-2x+a-1)\ph+(1-a)x+b-a$ be the $\ph$-expansion of $F(x)$. It follows that:
  \begin{enumerate}
		\item 
		If $a \not\equiv b\md4$ or  $a \equiv 3\md4$, then by Theorem \ref{intclos}, $2$ does not divide $i(K)$.
		\item
 If  $b \equiv a\equiv 1\md8$,  then $\nph{F}$ has  two sides joining $(0,v)$, $(1,1)$, and $(2,0)$ with $v\ge 3$. Thus $2\Z_K=\p_1\p_2$,  with  $\p_i$ a prime ideal of $\Z_K$ of residue degree $2$ for each $i=1,2$. 
 \item
 If  $b \equiv 1\md8$ and  $a\equiv 5\md8$,  then $\nph{F}$ has a single side joining $(0,2)$ and $(2,0)$, with residue degree $2$ and {$R_\la(F)=y^2+xy+x+1=(y-1)(y-x^2)$ in $\F_\ph[y]$. Thus $2\Z_K=\p_1\p_2$,  with $\p_i$ a prime of $\Z_K$ of residue degree $2$ for each $i$}. {In these cases, $\nu_2(i(K))=1$.}
		\item
 If  $b \equiv 5\md8$ and  $a\equiv 1\md8$,  then $\nph{F}$ has a single side joining $(0,2)$ and $(2,0)$, with residue degree $2$ and {$R_\la(F)=y^2+xy+1$, which is irreducible over $\F_\ph$.  Thus by Theorem \ref{ore}, $2\Z_K=\p^2$,  with $\p$ a prime ideal of $\Z_K$ of residue degree $4$}. Similarly, if  $b \equiv a\equiv 5\md8$,  then $\nph{F}$ has a single  side joining $(0,2)$, $(1,1)$, and $(2,0)$,  with residue degree $2$ and {$R_\la(F)=y^2+xy+x$, which is irreducible over $\F_\ph$. Thus  $2\Z_K=\p$,  with $\p$ a prime ideal of $\Z_K$ of residue degree $4$}.
 {In these cases, $\nu_2(i(K))=0$.}
 \end{enumerate}

  \end{enumerate}
$\hfill\Box$

\smallskip

\noindent
{\bf Proof of Theorem \ref{npib3}}.\\
 {By virtue of Lemma \ref{index}, we need to factorize $3\Z_K$ into powers of prime ideals of $\Z_K$. More precisely, $3$ divides $i(K)$ if and only if  $3\Z_K=\p_1\p_2\p_3\p_4$ with residue degree $1$ each prime ideal.}
Since $\triangle=2^4b(a^2-4b)^2$ is the discriminant of $F(x)$, if $3$ divides $(\Z_K:\Z[\al])$, then  $9$ divides $b$ or $3$ divides $a^2-4b$. If $3$  does not divide $b$, then $\ol{F(x)}$ has a square factor in $\F_3[x]$ if and only if $-a\equiv b\equiv 1 \md3$ or $a\equiv b\equiv 1 \md3$. Else $\ol{F(x)}$ is square free, and so by Dedekind's criterion,  $3$ does not divide $(\Z_K:\Z[\al])$. Recall that  $3$ is  a  divisor of $K$ if and only if there  are four prime ideals of $\Z_K$ lying above $3$.
\begin{enumerate}
		\item 
		If $-a\equiv b\equiv 1 \md3$, then  $ \ol{F(x)}= (x^2+1)^2$ in $\F_3[x]$. Thus there are at most two prime ideals of $\Z_K$ lying above $3$.
		\item 
	If $a\equiv b\equiv 1 \md3$, then $\nu_3(a^2-4b)\ge 1$.
		 Since $F'(x)=2x(2x^2+a)$  has three simple roots in $\F_3$, namely $0$, $1$ and $-1$, then by Hensel's lemma, let $s\in \Z$ such that $2s^2+a\equiv 0 \md{3^r}$ with $r=\nu_3(\triangle)+1$. Then  $ \ol{F(x)}= (x-s)^2(x+s)^2$ in $\F_3[x]$. Let $\ph_1=x-s$, $\ph_2=x+s$,  $F(x)=\ph_1^4+4s\ph_1^3+(6s^2+a)\ph_1^2+ (2as+4s^3)\ph_1+(b+as^2+s^4)$ and     $F(x)=\ph_2^4-4s\ph_2^3+(6s^2+a)\ph_2^2- (2as+4s^3)\ph_2+(b+as^2+s^4)$. 	 
						 Then $4(b+as^2+s^4)=2^2F(s)\equiv -(a^2-4b)\md{3^r}$ and $ F'(s)\equiv 0\md{3^r}$. So, $N^-_{\ph_i}(F)=S_i$ has a single side joining $(0, \nu_3(a^2-4b)$ and $(2,0)$. It follows that if $\nu_3(a^2-4b)$ is odd, then the degree of $S_i$ is $1$,  and so
						 $\ph_i$ provides two prime ideals of $\Z_K$ lying above $3$. If $\nu_3(a^2-4b)$ is even, then  $R_{\la_i}(F)(y)=y^2-(a^2-4b)_3$ for every $i=1,2$. Therefore, each $\ph_i$ provides two prime ideals of $\Z_K$ lying above $3$ if and only if  $(a^2-4b)_3\equiv 1\md3$.	
\item
		If $\nu_3(b)\ge 1$ and $a\equiv 1\md3$, then $ \ol{F(x)}= x^2(x^2+1)$ in $\F_3[x]$. Since $x^2+1$ is irreducible  in $\F_3[x]$, then there are at most three prime  ideals of $\Z_K$ lying above $3$. 
		\item
		If $\nu_3(b)= 2k+1$ is odd and $a\equiv -1\md3$, then $ \ol{F(x)}= x^2(x+1)(x-1)$ in $\F_3[x]$. Since $\nu_3(b)=1$, then the factor $x$ provides a unique prime  ideal of $\Z_K$ lying above $3$.  Thus,  there are at most three prime  ideals of $\Z_K$ lying above $3$. 
		\item
		If $\nu_3(b)= 2k$ for an odd natural integer $k$ and $a\equiv -1\md3$, then $ \ol{F(x)}= x^2(x+1)(x-1)$ in $\F_3[x]$. Since $\nu_3(b)=2k$, $N_x^-(F)$ has a single side  with $R_\la(F)(y)=y^2+b_3$ the attached residual polynomial of $F(x)$. Thus the factor $x$ provides two  prime  ideals of $\Z_K$ lying above $3$ if and only if $b_3\equiv -1\md3$. In this case there are  four prime  ideals of $\Z_K$ lying above $3$, and so $3$ divides $i(K)$. 
		\item
		Now assume that $9$ divides $b$ and $\nu_3(a)\ge 1$. In this case  $ \ol{F(x)}= x^4$ in $\F_3[x]$. It follows that:\\ 
		If   $2\nu_3(a)<\nu_3(b)$, then $\nu_3(a)=1$, $\nu_3(b)\ge 3$, and $\npp{F}=S_1+S_2$ has two sides, where $\ph=x$. Let $d_i$ be the degree of $S_i$ for every $i=1,2$. Then $d_2=1$ and $d_1\in\{1,2\}$. Thus there are at most three prime ideals of $\Z_K$ lying above $3$.\\
		If $2\nu_3(a)>\nu_3(b)$, then $\nu_3(b)\in \{2,3\}$ and  $\nu_3(a)\ge 2$, and $\npp{F}=S_1$ has a single side, where $\ph=x$. Let $d$ be the degree of $S_1$. Then  $d\in\{1,2\}$, and so there are at most two prime ideals of $\Z_K$ lying above $3$.
	\end{enumerate}
$\hfill\Box$

\vspace{0.5cm}

\noindent
{\bf Proof of Theorem  \ref{th1}}\\
We shall apply the main result of \cite{gppsys} (see also \cite{G19}) allowing to reduce the
index form equation in quartic number fields to a cubic form equation and to a pair of
quadratic form equations. 

Let $F(x)=x^4+ax^2+b$ an irreducible monogenic trinomial, then $(1,\alpha,\alpha^2,\alpha^3)$ is an 
integral basis in $K=\Q(\alpha)$, $\alpha$ being a root of $F(x)$. We can represent any $\gamma\in\Z_K$
in the form
\begin{equation}
\gamma=a+x\alpha+y\alpha^2+z\alpha^3
\label{gamma}
\end{equation}
with $a,x,y,z\in\Z$. We are going to determine all triples $(x,y,z)$, such that $\gamma$
generates a power integral basis in $K$ (for distinct $a\in\Z$ we obtain equivalent generators).
According to \cite{gppsys}, if $\gamma$ generates a power integral basis in $K$, then there exist 
$u,v\in \Z$ such that 
\begin{equation}
(u-av)(u^2-4bv^2)=\pm 1
\label{F}
\end{equation}
and
\begin{eqnarray}
Q_1(x,y,z)&=&x^2+ay^2-2axz+(a^2+b)z^2=u,\nonumber \\
Q_2(x,y,z)&=&y^2-xz+az^2 = v.  
\label{MN4}
\end{eqnarray}
By (\ref{F}) we have 
\[
u-av=\varepsilon,\;\;\; u^2-4bv^2=\delta
\]
with $\varepsilon=\pm 1,\delta=\pm 1$.
Then substituting $u=\varepsilon+av$ into $u^2-4bv^2=\delta$ we obtain
\begin{equation}
v^2(a^2-4b)+2\varepsilon av=\delta-\varepsilon^2=\delta-1.
\label{vv}
\end{equation}

\noindent
{\bf I.} If $\delta=-1$ then this implies that $v|(\delta-1)=-2$, hence 
$v=\pm 1$ or $v=\pm 2$.

\noindent
{\bf IA.} If $v=\pm 1$, then by $u=av+\varepsilon$ we have $u^2=a^2+2av\varepsilon+1$.
On the other hand we have $u^2-4bv^2=\delta =-1$, $u^2=4b-1$, whence
\[
a^2+2av\varepsilon+1=4b-1,
\]
\[
2av\varepsilon+2=4b-a^2,
\]
which can not be satisfied neither for even values of $a$ modulo 4, 
nor for odd values of $a$, modulo 2.

\noindent
{\bf IB.} If $v=\pm 2$, then by $u=av+\varepsilon$ we have $u^2=4a^2\pm 4a+1$.
On the other hand we have $u^2-4bv^2=\delta =-1$, $u^2=16b-1$, whence
\[
4a^2\pm 4a+1=16b-1,
\]
\[
4a^2\pm 4a-16b=-2,
\]
which can not be satisfied modulo 4.

\vspace{1cm}

\noindent
{\bf II.} If $\delta=1$ then (\ref{vv}) implies
\begin{equation}
v=0 \;\; {\rm or}\;\; v=\frac{-2\varepsilon a}{a^2-4b}.
\label{delta1}
\end{equation}
(Note that $a^2-4b\ne 0$, otherwise $K$ would only be a quadratic field.)\\

\noindent
{\bf IIA.}
In case 
\[
v=\frac{-2\varepsilon a}{a^2-4b}=(-\varepsilon)\frac{2a}{a^2-4b}
\]
we have
\[
u=av+\varepsilon=\frac{-2\varepsilon a^2}{a^2-4b}+\varepsilon
\]
whence
\[
u=(-\varepsilon)\frac{a^2+4b}{a^2-4b}.
\]
The above $u,v$ satisfies both $u-av=\varepsilon$ and $u^2-4bv^2=\delta=1$.
Set 
\[
u_0=\frac{u}{-\varepsilon}=\frac{a^2+4b}{a^2-4b},\;\; v_0=\frac{v}{-\varepsilon}=\frac{2a}{a^2-4b}
\]
($v\ne 0,v_0\ne 0$ since $a>1$).
We have
\[
a^2-4b=\frac{2a}{v_0}
\]
whence
\[
u_0=\frac{a^2+4b}{a^2-4b}=1+\frac{8b}{a^2-4b}=1+\frac{4bv_0}{a}
\]
($u_0\ne 1$ because $b>1$)
and 
\[
a=\frac{4bv_0}{u_0-1},\;\; a^2=\frac{16b^2v_0^2}{(u_0-1)^2}.
\]
Further, by $u_0=\frac{a^2+4b}{a^2-4b}$ we obtain
\[
u_0(a^2-4b)=a^2+4b,
\]
whence
\[
a^2(u_0-1)=4b(u_0+1),
\]
that is
\[
a^2=4b\frac{u_0+1}{u_0-1}.
\]
Comparing the two expessions obtained for $a^2$ we confer
\[
\frac{16b^2v_0^2}{(u_0-1)^2}=4b\frac{u_0+1}{u_0-1},
\]
whence
\[
b=\frac{u_0^2-1}{4v_0^2}
\]
and 
\[
a=\frac{4bv_0}{u_0-1}=\frac{u_0+1}{v_0}.
\]
All together, we obtain
\[
a=\frac{u\pm 1}{v},\;\; b=\frac{u^2-1}{4v^2}.
\]
Parameters of this type were excluded.\\

\noindent
{\bf IIB.}
Finally, if $v=0$ then equation (\ref{F}) implies $u=\pm 1$.
Observe that in our case
\[
Q_1(x,y,z)=x^2+ay^2-2xza+z^2(a^2+b)=(x-za)^2+ay^2+bz^2.
\]
Hence $Q_1(x,y,z)=u$ (see (\ref{MN4})) can only be satisfied for $u=1$ and then in view
of $a>1,b>1$ we have $x=\pm 1,y=0,z=0$, therefore up to equivalence $\alpha$ is the only generator of
power integral bases of $K$.
$\hfill\Box$

\vspace{1cm}

\noindent
{\bf Proof of Proposition \ref{ppprop}}\\
Let $\ph=x$. Then $\ol{F(x)}=\ph^4$ in $\F_2[x]$ and $\nph{F}=S$ has a single side of degree $\gcd(4,3)=1$. Thus  $F(x)$ is irreducible over $\Q_2$. Let $K$ be the number field generated by a root $\al$ of $F(x)$.  

Then there is a unique valuation $\om$ of $K$ extending $\nu_2$. 
Since $\nu_2(\Z_K:\Z[\al])=3$, we conclude that $F(x)$ is not a  monogenic polynomial. Now  let $\theta =\frac{\al^3}{4}$.
	Then $\theta \in K$. Since $3$ and $4$ are coprime, we conclude that $K=\Q(\theta)$.	 Let us show that $\Z_K=\Z[\theta]$, and so $K$ is monogenic.  By \cite[Corollary 3.1.4]{En}, in order to show that $\theta\in \Z_K$, we need to show that $\om(\theta)\ge 0$, where $\om$ is the unique valuation of $K$ extending $\nu_2$. Since $\nph{F}=S$ has  a single side of slope $-3/4$, we conclude that $\om(\al)=3/4$, and so $\om(\theta)=\frac{9}{4}-2=\frac{1}{4}$. Let $g(x)$ be the minimal polynomial of $\theta$ over $\Q$. By the formula relating roots and coefficients of a monic polynomial, we conclude that $g(x)=x^{4}+\sum_{i=1}^{4}(-1)^is_ix^{4-i}$, where $s_i=\displaystyle\sum_{k_1<\dots<k_i}\theta_{k_1}\cdots\theta_{k_i}$ and $\theta_{1},\dots, \theta_{4}$ are the $\Q_p$-conjugates of $\theta$. Since there is a unique valuation extending $\nu_2$ to any algebraic extension of $\Q_2$, we conclude that $\om(\theta_i)=1/4$ for every $i=1,\dots,4$. Thus $\nu_2(s_{4})=\om(\theta_{1}\cdots\theta_{4})=4\times 1/4=1$ and $\nu_2(s_{i})\ge i/4$ for every $i=1,\dots, 3$, which  means that $g(x)$ is a $2$-Eisenstein polynomial. Hence $2$ does not divide the index $(\Z_K:\Z[\theta])$. As $\tr(F)=2^{15}b(a^2-2b)^2$ and by definition of $\theta$, $2$ is the unique positive prime integer candidate to divide  $(\Z[\al]:\Z[\theta])$, we conclude that for every prime integer $p$, $p$ does not divide $(\Z_K:\Z[\theta])$, which means that $\Z_K=\Z[\theta]$.
$\hfill\Box$

\end{document}